\documentclass[11pt,reqno]{amsart}
\usepackage{amsmath,amsfonts, amssymb, amsthm}
  \usepackage{paralist}
  \usepackage{graphics} 
  \usepackage{epsfig} 
\usepackage{graphicx}  \usepackage{epstopdf}
 \usepackage[colorlinks=true]{hyperref}
\hypersetup{urlcolor=blue, citecolor=red}

\newtheorem{theorem}{Theorem}[section]

\newtheorem{lemma}[theorem]{Lemma}

\newtheorem{conjecture}{Conjecture}[section]

\theoremstyle{definition}

\newtheorem{remark}{Remark}[section]

\def\pmod #1{\ ({\rm{mod}}\ #1)}
\def\Z{\Bbb Z}
\def\N{\Bbb N}

\def\Q{\Bbb Q}

\def\l{\left}
\def\r{\right}
\def\bg{\bigg}
\def\({\bg(}
\def\){\bg)}
\def\t{\text}
\def\f{\frac}
\def\mo{{\rm{mod}\ }}
\def\pmod#1{\ (\mo\ #1)}

\def\cs{\ldots}

\def\ls{\leqslant}
\def\gs{\geqslant}
\def\se {\subseteq}
\def\sm{\setminus}
\def\bi{\binom}

\def\ve{\varepsilon}

\def\eq{\equiv}

\def\da{\delta}

\def\Proof{\noindent{\it Proof}}
\def\Ack{\medskip\noindent {\bf Acknowledgment}}

\begin{document}
\hbox{Submitted to J. Syst. Sci. Complex.}
\medskip

\title[Mixed quantifier prefixes over Diophantine equations]
      {Mixed quantifier prefixes over Diophantine equations with integer variables}
\author[Zhi-Wei Sun]{Zhi-Wei Sun}


\address{Department of Mathematics, Nanjing
University, Nanjing 210093, People's Republic of China}
\email{zwsun@nju.edu.cn}

\keywords{Undecidability, Diophantine equations, Hilbert's tenth problem, mixed quantifiers.
\newline \indent The research is supported by the National Natural Science Foundation of China (grant no. 12371004), and the main results date back to the author's PhD thesis (Nanjing University, 1992).}

\begin{abstract} In this paper we first review the history of Hilbert's Tenth Problem, and then study mixed quantifier prefixes over Diophantine equations with integer variables. For example, we prove that $\forall^2\exists^4$ over $\mathbb Z$ is undecidable, that is, there is no algorithm to determine for any $P(x_1,\ldots,x_6)\in\mathbb Z[x_1,\ldots,x_6]$ whether
$$\forall x_1\forall x_2\exists x_3\exists x_4\exists x_5\exists x_6(P(x_1,\ldots,x_6)=0),$$
where $x_1,\ldots,x_6$ are integer variables.
We also have some similar undecidable results with universal quantifies bounded, for example,
$\exists^2\forall^2\exists^2$ over $\mathbb Z$ with $\forall$ bounded is undecidable.
We conjecture that $\forall^2\exists^2$ over $\mathbb Z$ is undecidable.
\end{abstract}
\maketitle

\section{Introduction}

In 1900, at the Paris conference of ICM, D. Hilbert presented 23 famous mathematical problems.
Many of them are questions of others, however the tenth one is due to Hilbert himself.
In modern language, {\it Hilbert's Tenth Problem} (HTP) asks for an effective algorithm to test whether
an arbitrary polynomial equation
$$P(z_1,\cs,z_n)=0$$
(with integer coefficients) has solutions over the ring $\Z$ of the integers, where $n$ is an arbitrary positive integer.
However, the concept of algorithm or computation was vague in 1900.

Let $\N=\{0,1,2,\ldots\}$ and call each $n\in\N$ a {\it natural number}.
What kind of number-theoretic functions into $\N$ (with natural number variables)
are computable? This was investigated by logicians in the 1930s.

We first introduce the basic functions:

(1) {\it Zero function}: $O(x)=0$ (for all $x\in\N$).

(2) {\it Successor function}: $S(x)=x+1$.

(3) {\it Projection functions}: $I_{nk}(x_1,\ldots,x_n)=x_k\ (1\ls k\ls n)$

For number-theoretic functions $g(y_1,\ldots,y_m)$ and $h_i(x_1,\ldots,x_n)\ (1\ls i\ls m)$, we define their {\it composition} as follows:
$$f(x_1,\ldots,x_n)=g(h_1(x_1,\ldots,x_n),\ldots,h_m(x_1,\ldots,x_n))$$

Given number-theoretic functions $g(x_1,\ldots,x_n)$
and $h(x_1,\ldots,x_n,y,z)$, we define
$$\begin{cases}f(x_1,\ldots,x_n,0)=g(x_1,\ldots,x_n),
\\f(x_1,\ldots,x_n,y+1)=h(x_1,\ldots,x_n,y,f(x_1,\ldots,x_n,y)),
\end{cases}$$
and say that $f$ is obtained from $g$ and $h$ via {\it primitive recursion}.

For a number-theoretic function $g(x_1,\ldots,x_n,y)$, we define
$$f(x_1,\ldots,x_n)=\mu y\gs0\,(g(x_1,\ldots,x_n,y)=0)$$
as the least $y\in\N$ with $g(x_1,\ldots,x_n,y)=0$;
if  $g(x_1,\ldots,x_n,y)\not=0$ for all $y\in\N$, then $f(x_1,\ldots,x_n)$ is undefined.
We say that $f$ is obtained from the function $g$ via the {\it $\mu$-operator}.

 The {\it partial recursive functions} are the basic functions and those obtained from the basic functions by applying composition, primitive recursion and the $\mu$-operator a finite number of times.
For any partial recursive function $f$, it is easy to see that if $f(x_1,\ldots,x_n)$ is defined then the value $f(x_1,\ldots,x_n)$ is effectively computable by intuition.

In 1936 A. Turing introduced the notion of Turing machine which is an abstract machine that manipulates symbols on a infinite strip of tape according to a finite table of rules (i.e., a program) involving four kinds of basic operations: Write $1$, change $1$ to $0$ (blank), move to the left unit (L), move to the right unit (R).
A function $f(x_1,\ldots,x_n)$ is {\it Turing computable} if there is a program according to which the Turing machine with initial inputs $x_1,\ldots,x_n$
finally stops and yields the value $f(x_1,\ldots,x_n)$ as output if $f(x_1,\ldots,x_n)$ is defined, and never stops if $f(x_1,\ldots,x_n)$ is undefined.

The partial recursive functions and Turing computable functions were proved to be equivalent by S. C. Kleene in 1936. The following thesis was proposed by A. Church in the same year.
\smallskip

{\bf Church's Thesis}. {\it If a function $f$ into $\N$ with natural number variables is effectively computable by intuition, then it must be a partial recursive function (equivalently, a Turing computable function)}.
\smallskip

Church's Thesis has been widely accepted after 1936. So we have the exact definition of computable functions which refer to the partial recursive functions or Turing computable functions, and hence
HTP has its accurate meaning.

A subset $A$ of $\N$ is said to be an {\it r.e. (recursively enumerable) set} (or {\it a semi-decidable set}) if the function
$$f_A(x)=\begin{cases}1&\mbox{if}\ x\in A,
\\\t{undefined}&\t{if}\ x\in\N\sm A.
\end{cases}$$
is a partial recursive function. It is easy to show that
\begin{align*}&A\se\N\ \text{is an r.e. set}
 \\\iff&A\ \t{is the domain of a partial recursive function}
 \\\iff&A\ \t{is the emptyset or the range of a total recursive function}\ f(x).
 \end{align*}

 A set $A\se \N$ is called {\it decidable} or {\it recursive}, if
the characteristic function
$$\chi_A(x)=\begin{cases}1&\t{if}\ x\in A,
\\0&\t{if}\ x\in\N\sm A.\end{cases}$$
is Turing computable (or recursive).
Clearly, $A$ is recursive if and only if both $A$ and $\N\setminus A$ are r.e. sets.
A well known result in the theory of computability states that there is a nonrecursive r.e. set
(cf. \cite[pp.\,140-141]{C80}).

From now on, variables range over $\Z$ unless specified.
Let $P(z_1,\ldots,z_n)\in\Z[z_1,\ldots,z_n]$. Then
\begin{align*}&\exists z_1\ldots\exists z_n\,(P(z_1,\ldots,z_n)=0)
\\\iff&\exists x_1\gs0\ldots \exists x_n\gs0\,\bigg(\prod_{\ve_1,\ldots,\ve_n\in\{\pm1\}}P(\ve_1x_1,\ldots,\ve_nx_n)=0\bigg).
\end{align*}
On the other hand, by Lagrange's four-square theorem (each $m\in\N$ can be written as the sum of four squares), we have
\begin{align*} &\exists x_1\gs0\ldots\exists x_n\gs0\,(P(x_1,\ldots,x_n)=0)
\\\iff&\exists u_1\exists v_1\exists y_1\exists z_1\ldots\exists u_n\exists v_n\exists y_n\exists z_n
\\&(P(u_1^2+v_1^2+y_1^2+z_1^2,\ldots,u_n^2+v_n^2+y_n^2+z_n^2)=0)
\end{align*}
So HTP has the following equivalent form (HTP over $\N$): Is there an algorithm to decide for any polynomial $P(x_1,\ldots,x_n)$
 with integer coefficients whether the Diophantine equation
 $P(x_1,\ldots,x_n)=0$
 has solutions with $x_1,\ldots,x_n\in\N$?

 A relation $R(a_1,\ldots,a_m)$ with $a_1,\ldots,a_m\in\N$ is said to be {\it Diophantine} if
there is a polynomial $P(t_1,\ldots,t_m,x_1,\ldots,x_n)$ with integer coefficients
such that
$$R(a_1,\ldots,a_m)\iff\exists x_1\gs0\ldots\exists x_n\gs0\,(P(a_1,\ldots,a_m,x_1,\ldots,x_n)=0).$$
A set $A\se\N$ is Diophantine if and only if the predicate $a\in A$ is Diophantine.
It is easy to see that
any Diophantine set $A$ is an r.e. set. In fact, for a given element $a\in \N$
we may search for the natural number solutions of the Diophantine equation associated with $A$.
If it has a solution, then we will find one and let the computer stop and give the output $1$. If it has no solution, the computer will never stop.

 In 1944 E. L. Post  thought that HTP begs for an unsolvability proof, i.e., HTP might be undecidable
(cf. \cite{D73,M-book}).
In 1953, M. Davis \cite{D53} published the following hypothesis.
\smallskip

{\bf Davis Daring Hypothesis}. {\it Any r.e. set $A\se\N$ is Diophantine}.
\smallskip

Under this hypothesis, for any nonrecursive r.e. set $A$ there is a polynomial $P(x,x_1,\ldots,x_n)$
such that for any $a\in\N$ we have
$$a\in A\iff \exists x_1\gs0\ldots\exists x_n\gs0\,(P(a,x_1,\ldots,x_n)=0).$$
Thus Davis Daring Hypothesis implies that HTP over $\N$ is undecidable.

The {\it exponential Diophantine equations} over $\N$ have the form
$$E_1(x_1,\ldots,x_m)=E_2(x_1,\ldots,x_m),$$
where $E_1$ and $E_2$ are expressions constructed from variables and particular natural numbers using
addition, multiplication, and exponentiation. Here is an example of exponential Diophantine equation:
$$x^{2^y}+y^2+y^{y^z}=5z^{x^x+3z}.$$
A relation $R(a_1,\ldots,a_m)$ with $a_1,\ldots,a_m\in\N$ is said to be {\it exponential Diophantine} if
there is an exponential Diophantine equation
 $$E(t_1,\ldots,t_m,x_1,\ldots,x_n)=0$$
over $\N$ such that
$$R(a_1,\ldots,a_m)\iff \exists x_1\gs0\ldots\exists x_n\gs0\,(E(a_1,\ldots,a_m,x_1,\ldots,x_n)=0).$$
A set $A\se\N$ is called exponential Diophantine if the predicate $a\in A$ is Diophantine.
The following important result concerning exponential Diophantine equations was established by Davis, H. Putnam and J. Robinson \cite{DPR} in 1961.
\smallskip

{\bf Davis-Putnam-Robinson Theorem}.  {\it Any r.e. set is exponential Diophantine.
Thus there is no algorithm to decide for any given exponential Diophantine equation whether it has solutions over $\N$}.
\smallskip

Based on this result, in 1970 Y. Matiyasevich \cite{M70} utilized the Fibonacci sequence to prove that the exponential relation $a=b^c$ (with $a,b,c\in\N$) is Diophantine. This, together with the Davis-Putnam-Robinson Theorem, led him to prove the Davis Daring Hypothesis completely. Thus, HTP was finally solved negatively in 1970.
The reader may consult Davis \cite{D73} for a popular introduction to the negative solution of HTP.

In 1975 Matiyasevich proved further that
 any r.e. set $A\se\N$ has a Diophantine representation over $\N$ with only 9 unknowns,
 the detailed proof of this 9 unknowns theorem appeared in J. P. Jones \cite{J82}.

 Note that a system of finitely many Diophantine equations over $S\se\Z$ is equivalent to a single Diophantine equation over $S$.
In fact, if $P_i(z_1,\ldots,z_n)\in\Z[z_1,\ldots,z_n]$ for all $i=1,\ldots,k$, then
\begin{align*}&P_1(z_1,\ldots,z_n)=0\land \ldots\land P_k(z_1,\ldots,z_n)=0
\\\iff&P_1^2(z_1,\ldots,z_n)+\ldots+P_k^2(z_1,\ldots,z_n)=0.\end{align*}

 For $i=1,\ldots,n$, let each $\rho_i$ be one of the two quantifiers $\forall$ and $\exists$.
 If there is no algorithm to determine for any $P(x_1,\ldots,x_n)\in\Z[x_1,\ldots,x_n]$
 whether  $$\rho_1 x_1\gs0\cdots\rho_n x_n\gs0\,(P(a,x_1,\ldots,x_n)=0),$$
 then we say that $\rho_1\cdots\rho_n$ over $\N$ is undecidable.
 For example, $\exists^9$ over $\N$ is undecidable by the 9 unknowns theorem, but it is open whether $\exists^8$ over $\N$ is decidable or not.
 We may also consider $\rho_1\cdots\rho_n$ over $\N$ with $\forall$ bounded, for example, Matiyasevich \cite{M76} proved that $\exists\forall\exists^2$ over $\N$ with $\forall$ bounded is undecidable, that is, there is no algorithm to determine for any  $P(x)\in \Z[x]$ and $Q(x_1,\ldots,x_4)\in\Z[x_0,\ldots,x_4]$ whether
 $$\exists x_1\gs0\forall x_2\in[0,P(x_1)]\exists x_3\gs0\exists x_4\gs0\,(Q(a,x_1,\ldots,x_4)=0).$$
 After the negation solution of Hilbert's tenth problem, it is natural to ask the following question: For what kinds of mixed quantifier prefixes $\rho_1\cdots\rho_n$, $\rho_1,\cdots\rho_n$ over $\N$
(with $\forall$ bounded or unbounded) is undecidable?
After a series of efforts due to Matiyasevich \cite{M73,M76}, Matiyasevich and Robinson \cite{MR74,MR},
and Jones \cite{J81}, the only open cases are
$\forall\exists^2$, $\exists\forall\exists$, and $\exists\forall\exists$ with $\forall$ bounded.
J. M. Rojas \cite[Conjecture 3]{Rojas} conjectured that $\exists\forall\exists$ over $\N$ is decidable.

Both $\exists$ over $\N$ and $\exists$ over $\Z$ are decidable in polynomial time (see, e.g., \cite[p.\,525]{MR}).
In fact, if $a_0,a_1,\ldots,a_n$ and $z$ are integers
with $a_0z\not=0$ and $\sum_{i=0}^na_iz^{n-i}=0$, then
$$|z|^n\ls|a_0z^n|\ls\sum_{i=1}^n|a_i|\cdot|z|^{n-i}\ls \sum_{i=1}^n|a_i|\cdot|z|^{n-1}$$
and hence
$$|z|\ls \sum_{i=1}^n|a_i|.$$

In 1987 S. P. Tung \cite{T87} showed that for each $n\in\Z^+=\{1,2,3,\ldots\}$
the problem to determine
$$\forall x_1\cdots\forall x_n \exists x_{n+1}\,(P(x_1,\ldots,x_n,x_{n+1})=0)$$
with $P$ a general polynomial in $\Z[x_1,\ldots,x_{n+1}]$ is co-NP-complete.

For a finite sequence of quantifiers $\rho_1,\ldots,\rho_n$, we say that $\rho_1\cdots\rho_n$
over $\Z$ is undecidable if there is no algorithm to determine for any $P(x_1,\ldots,x_n)\in\Z[x_1,\ldots,x_n]$ whether
\begin{equation}\label{rho}\rho_1x_1\cdots\rho_nx_n\,(P(x_1,\ldots,x_n)=0).
\end{equation}
What kinds of $\rho_1\cdots\rho_n$ over $\Z$
are undecidable? In 1985 Tung \cite{T85} proved that $\exists^{27}$ and $\forall^{27}\exists^2$ over $\Z$
are undecidable.
We may also consider $\rho_1\cdots\rho_n$
over $\Z$ with $\forall$ bounded. We say that $\rho_1\cdots\rho_n$
over $\Z$ with $\forall$ bounded is undecidable if there is no general algorithm to determine
whether \eqref{rho} with $\rho_j x_j$ (for those $1\ls j\ls n$ with $\rho_j=\forall$)
replaced by
$$\forall x_j\in[P_j(x_i:1\ls i<j\ \&\ \rho_i=\exists),Q_j(x_i:1\ls i<j\ \&\ \rho_i=\exists)]$$
holds or not, where $P$ and those $P_j$ and $Q_j$ with $\rho_j=\forall$ are polynomials with integer coefficients.  For example, $\exists\forall^2\exists$ over $\Z$ is undecidable if and only if
there is no algorithm to determine for $P_1(x),P_2(x),P_3(x),P_4(x)\in\Z[x]$ and $Q(x_1,\ldots,x_4)\in\Z[x_1,\ldots,x_4]$ whether
\begin{align*}\exists x_1\forall x_2\in[P_1(x_1),P_2(x_1)]\forall x_3\in[P_3(x_1),P_4(x_1)]\exists x_4
(Q(x_1,x_2,x_3,x_4)=0).
\end{align*}
Clearly, if $\rho_1\cdots \rho_n$ over $\Z$ (with $\forall$ bounded or not) is decidable, then so is
$\rho_{i_1}\rho_{i_2}\cdots\rho_{i_m}$ over $\Z$ with $1\ls i_1<i_2<\ldots<i_m$.
Also, $\forall \rho_1\cdots \rho_n$ over $\Z$ with $\forall$ bounded is decidable if and only if
$\rho_1\cdots \rho_n$ over $\Z$ with $\forall$ bounded is decidable.
Note also that
$\rho_1,\cdots\rho_n$ over $\Z$ is undecidable (with $\forall$ unbounded or bounded) if and only if
$\rho_1\cdots\rho_n\forall$ over $\Z$ is undecidable. In fact, for the polynomial
$$P(x_1,\ldots,x_n,t)=\sum_{i=0}^kP_i(x_1,\ldots,x_n)t^i\ \t{with}\ P_i(x_1,\ldots,x_n)\in\Z[x_1,\ldots,x_n],$$
we have
$$\forall t(P(x_1,\ldots,x_n,t)=0)\iff \sum_{i=0}^k P_i(x_1,\ldots,x_n)^2=0$$
for all $x_1,\ldots,x_n\in\Z$; if $P_*(x_1,\ldots,x_n),P^*(x_1,\ldots,x_n)\in\Z[x_1,\ldots,x_n]$,
then for any $x_1,\ldots,x_n\in\Z$ we have
\begin{align*}&\ \forall t\in[P_*(x_1,\ldots,x_n),P^*(x_1,\ldots,x_n)](P(x_1,\ldots,x_n,t)=0)
\\\iff&\ \sum_{0\ls r\ls P^*(x_1,\ldots,x_n)-P_*(x_1,\ldots,x_n)}\(\sum_{i=0}^kP_i(x_1,\ldots,x_n)(P_*(x_1,\ldots,x_n)+r)^i\)^2=0
\\\iff&\ \sum_{0\ls r\ls P^*(x_1,\ldots,x_n)-P_*(x_1,\ldots,x_n)}\sum_{j=0}^{2k}\bar P_j(x_1,\ldots,x_n)r^j=0,
\end{align*}
where $\bar P_j(x_1,\ldots,x_n)\ (0\ls j\ls 2k)$ are suitable polynomials with integer coefficients.
For any $j,m\in\N$, it is well known (cf. \cite[pp.\,228-231]{IR}) that
$$\sum_{r=0}^m r^j=\f1{j+1}\sum_{r=0}^m\l(B_{j+1}(r+1)-B_{j+1}(r)\r)=\f{B_{j+1}(m+1)-B_{j+1}(0)}{j+1},$$
where $B_{j+1}(x)$ denotes the Bernoulli polynomial (with rational number coefficients) of degree $j+1$.

Based on the work \cite{S92a}, the author \cite{S21} proved that for any r.e. set $A$ there is a polynomial
$P(x_0,\ldots,x_9)\in\Z[x_0,\ldots,x_9]$ such that for any $a\gs0$ we have
$$a\in A\iff \exists x_1\cdots\exists x_8\exists x_9\gs0\,(P(a,x_1,\ldots,x_9)=0).$$
(See also the book \cite{S24} for a complete proof.)
This implies Matiyasevich's 9 unknowns theorem since
$$a\in A\iff \exists x_1\gs0\cdots\exists x_9\gs0\,\(\prod_{\ve_1,\ldots,\ve_8\in\{\pm1\}}
P(a,\ve_1x_1,\ldots,\ve_8x_8,x_9)=0\).$$
As a consequence
of this result, the author \cite{S21} obtained  the 11 unknowns theorem ($\exists^{11}$
over $\Z$ is undecidable) and also the undecidability of $\forall^{10}\exists^2$
and $\forall^9\exists^3$ over $\Z$. For their applications, one may consult
\cite{Da,MPR,MS,RGK,ZS}.
The author \cite[Conjecture 1.8]{S21} conjectured that there is no algorithm to determine for any $P(x,y,z)\in\Z[x,y,z]$ whether the equation $P(x^2,y^2,z^2)=0$ has integer solutions, which implies that $\exists^3$ over $\Z$ is undecidable.

Now we state the main result of this paper.

\begin{theorem} \label{Th1.0} {\rm (i)} All those
\begin{gather*}\forall\exists^7,\ \forall^2\exists^4,\ \exists\forall\exists^4,\ \exists\forall^2\exists^3,\ \exists^2\forall\exists^3,\ \forall\exists\forall\exists^3,
\\ \forall\exists^2\forall^2\exists^2,\ \forall^2\exists\forall^2\exists^2,
\forall\exists\forall^3\exists^2,\ \exists^2\forall^3\exists^2,\ \exists\forall\exists\forall^2\exists^2,\ \exists\forall^6\exists^2
\end{gather*}
over $\Z$ are undecidable.

{\rm (ii)} All those
$$\exists\forall\exists^4,\ \exists\forall^2\exists^3,\ \exists^2\forall\exists^3,
\ \exists^2\forall^2\exists^2, \exists\forall\exists\forall\exists^2,\ \exists\forall^5\exists^2$$
over $\Z$ with $\forall$ bounded are undecidable.
\end{theorem}

\begin{remark} In 1991 the author learnt from Tung that R. M. Robinson was the first person
to ask for such undecidable results over $\Z$.
\end{remark}

Given a finite sequence of quantifiers $\rho_1,\ldots,\rho_n$, we say that a set $A\se\N$
has a $\rho_1\cdots\rho_n$-representation over $\Z$ if there is a polynomial $P(x_0,\ldots,x_n)
\in\Z[x_0,\ldots,x_n]$ such that for any $a\in\N$ we have
$$a\in A\iff \rho_1x_1\cdots\rho_n x_n\,(P(a,x_1,\ldots,x_n)=0).$$
Similarly, we may define $\rho_1\cdots\rho_n$-representations over $\Z$ with $\forall$ bounded.
The author \cite{S21} actually proved that any r.e. set $A\se \N$
has a $\exists^{11}$-representation over $\Z$, and any
co-r.e. set (i.e.,
the complement of an r.e. set $A\se\N$) has a  $\forall^{10}\exists^2$-representation,
and a $\forall^9\exists^3$-representation over $\Z$.
By B.-K. Oh and the author \cite[Corollary 1.1]{OS}, the set
$$S=\{2n+1:\ n\in\Z^+,\ \t{and}\ 2n+1\ \t{is not a prime congruent to}\ 3\ \t{modulo}\ 4\}$$
has a surprising $\exists^3$-representation over $\Z$: $a\in\N$ belongs to $S$ if and only if
$$\exists x\exists y\exists z\,(a^2=(2x+1)^2+8(2y+1)^2+8(2z+1)^2).$$

For a subset $A$ of $\N$, we write $\bar A$ for $\N\sm A$, the complement of $A$ in $\N$.
Theorem \ref{Th1.0} follows immediately from our following three
theorems.

\begin{theorem}\label{Th1.1} Let $A\se\N$ be an r.e. set.

{\rm (i)} $A$ has a $\exists^2\forall\exists^3$-representation over $\Z$. Also, we may replace
$\exists^2\forall\exists^3$ by either of
 $\exists^2\forall^3\exists^2$ and $\exists\forall\exists\forall^2\exists^2. $
Also, $\bar A$ has a $\forall^2\exists^4$-representation, a $\forall\exists\forall\exists^3$-representation, and a
$\forall^2\exists\forall^2\exists^2$-representation
over $\Z$.

{\rm (ii)} $A$ has a $\exists^2\forall^2\exists^2$-representation over $\Z$ with $\forall$ bounded.
Also, we may replace
$\exists^2\forall^2\exists^2$ by either of
$\exists^2\forall \exists^3$ and $\exists\forall\exists\forall\exists^2.$
\end{theorem}

\begin{theorem}\label{Th1.2} Let $A\se\N$ be an r.e. set.

{\rm (i)} $A$ has a $\exists\forall\exists^4$-representation over $\Z$, and also a
$\exists\forall\exists^4$-representation over $\Z$
with $\forall$ bounded.

{\rm (ii)} $\bar A$ has a $\forall\exists\forall^3\exists^2$-representation over $\Z$.
\end{theorem}

\begin{theorem}\label{Th1.3} Let $A\se\N$ be an r.e. set.

{\rm (i)} $A$ has a $\exists\forall^2\exists^3$-representation over $\Z$,
and also a $\exists\forall^2\exists^3$-representation over $\Z$ with $\forall$ bounded.
Also, $A$ has a $\exists\forall^6\exists^2$-representation over $\Z$,
and a $\exists\forall^5\exists^2$-representation over $\Z$ with $\forall$ bounded.

{\rm (ii)} $\bar A$ has a $\forall \exists^7$-representation over $\Z$ and also a
$\forall\exists^2\forall^2\exists^2$-representation over $\Z$.
\end{theorem}

In Section 2 we will prove an auxiliary theorem. Sections 3-5 will be devoted to our proofs
of Theorems \ref{Th1.1}-\ref{Th1.3}, respectively.

Although we have Theorem \ref{Th1.0}, there are many finite sequences
of quantifiers (such as $\exists\forall^m\exists\forall^n\exists$ ($m\in\{2,3,4\}$ and $n\in\N$) with $\forall$ bounded or not) for which we don't know whether they are undecidable over $\Z$.

It is believed that $\exists^3$ over $\Z$ might be undecidable (cf. \cite{MR}). We pose here a conjecture for further research.

\begin{conjecture} $\forall^2\exists^2$ over $\Z$ is undecidable.
\end{conjecture}

We mention that HTP over the field $\Q$ of rational numbers is a difficult open problem.
Also, for a general number field $K$ (which is a finite extension of the field $\Q$), HTP over
the ring $O_K$ of all algebraic integers in $K$, remains open. The reader may consult \cite{Da,Sh} for certain progress.

Our theorems in this paper should have some potential applications. For example, if one investigates
mixed quantifiers over Diophantine equations with variables ranging over the rational field $\Q$ or the Gaussian ring $\Z[i]$, our results on mixed quantifiers over $\Z$ will be useful.

\section{An auxiliary theorem}
 \setcounter{equation}{0}
 \setcounter{conjecture}{0}
 \setcounter{theorem}{0}

In this section we adapt Matiyasevich and Robinson's ideas in \cite{MR74,MR} to establish an auxiliary theorem on representations of r.e. sets over $\Z$
which will be helpful to our later proofs of Theorems \ref{Th1.1}-\ref{Th1.3}.

 \begin{lemma} \label{Lem2.1} Let $B\gs b>0$ and $0<n_0<\ldots<n_{\nu}$. Then an integer $c$ has the form $\sum_{i=0}^{\nu}z_iB^{n_i}$ with $z_i\in\{0,\ldots,b-1\}$ for all $i=0,\ldots,\nu$ if and only if
 every interval $[\sigma_i,\tau_i]\ (i=0,\ldots,\nu+1)$ contains at least an integer, where
 \begin{gather*}\sigma_0=\f c{B^{n_0}},\ \sigma_i=\f{c+1-bB^{n_{i-1}}}{B^{n_i}}\ (i=1,\ldots,\nu),
 \ \sigma_{\nu+1}=\f{c+1-bB^{n_{\nu}}}{(b^2+c^2)B^{n_{\nu}}},
 \\\tau_i=\f c{B^{n_i}}\ (i=0,1,\ldots,\nu)\ \t{and}\ \tau_{\nu+1}=\f c{(b^2+c^2)B^{n_{\nu}}}.
 \end{gather*}
 \end{lemma}
 \Proof. We first prove the ``only if" direction. Suppose that $c=\sum_{i=0}^{\nu}z_iB^{n_i}$ with $z_i\in\{0,\ldots,b-1\}$ for all $i=0,\ldots,\nu$.
 For any $j=0,\ldots,\nu$, we have
 \begin{align*}\sum_{i=0}^jz_iB^{n_i}\ls& (b-1)B^{n_j}+\sum_{0\ls i<j}(b-1)B^{n_i}
 \\\ls& (b-1)B^{n_j}+\sum_{k=0}^{n_j-1}(B-1)B^k=(b-1)B^{n_j}+B^{n_j}-1=bB^{n_j}-1
 \end{align*}
 In particular, $0\ls c=\sum_{i=0}^{\nu}z_iB^{n_i}\ls bB^{n_{\nu}}-1$ and hence $\sigma_{\nu+1}\ls 0\ls\tau_{\nu+1}$. Set
 $$x_i=\sum_{j=i}^{\nu}z_jB^{n_j-n_i}\quad\ \ \t{for}\ i=0,1,\ldots,\nu.$$
 Then $x_0=\sigma_0$,  and
 $$x_i=\f{c-\sum_{j=0}^{i-1}z_jB^{n_j}}{B^{n_i}}\gs\f{c-(bB^{n_{i-1}}-1)}{B^{n_i}}=\sigma_i$$
 for all $i=1,\ldots,\nu$. Also,
 $$x_i\ls\sum_{j=0}^\nu z_jB^{n_j-n_i}=\f c{B^{n_i}}=\tau_i$$
 for all $i=0,\ldots,\nu$. Therefore, each interval $[\sigma_i,\tau_i]\ (i=0,\ldots,\nu+1)$
 contains an integer. This proves the ``only if" direction.

 Now we consider the ``if" direction. Suppose that there are integers $x_0,\ldots,x_{\nu+1}$ with $\sigma_i\ls x_i\ls\tau_i$ for all $i=0,\ldots,\nu+1$. Since
 $$|c+1-bB^{n_{\nu}}|\ls bB^{n_{\nu}}-1+|c|<(b^2+c^2)B^{n_{\nu}},$$
 we have $|\sigma_{\nu+1}|<1$. Note also that $|\tau_{\nu+1}|<1$. As $$-1<\sigma_{\nu+1}\ls x_{\nu+1}\ls\tau_{\nu+1}<1,$$ we must have $x_{\nu+1}=0$. From $\sigma_{\nu+1}\ls0\ls\tau_{\nu+1}$,
 we get $0\ls c<bB^{n_\nu}\ls B^{n_\nu+1}$. No matter $B>1$ or $B=1$, we can write
 $$c=\sum_{k=0}^{n_\nu}c_kB^k$$
 with $c_{n_\nu}\in\{0,\ldots,b-1\}$ and $c_k\in\{0,\ldots,B-1\}$ for all $k=0,\ldots,n_{\nu}-1$.

 As $\sigma_0=\tau_0$, we have $\sigma_0=x_0\in\Z$ and hence $B^{n_0}\mid c$.
 Thus $c_k=0$ for all $k=0,\ldots,n_{0}-1$.

 Let $1\ls i\ls\nu$. As $\sigma_i\ls x_i\ls\tau_i$, we have
 $$0\ls c-x_iB^{n_i}\ls bB^{n_{i-1}}-1\ls BB^{n_{i-1}}-1<B^{n_{i-1}+1}\ls B^{n_i}$$
 and hence $c-x_iB^{n_i}$ is the least nonnegative residue of $c$ modulo $B^{n_i}$.
 Thus
 $$\sum_{k=0}^{n_i-1}c_kB^k=c-x_iB^{n_i}< bB^{n_{i-1}}.$$
 It follows that $c_{n_{i-1}}<b$ and $c_{n_{i-1}+1}=\ldots=c_{n_i-1}=0$.

 By the above, we have $c=\sum_{i=0}^{\nu} z_iB^{n_i}$ with $z_i=c_{n_i}\in\{0,\ldots,b-1\}$
 for all $i=0,\ldots,\nu$. This ends our proof of the ``if" direction. \qed

\begin{lemma}\label{Lem2.2} Let $\sigma_0,\tau_0,\ldots,\sigma_k,\tau_k$ be real numbers
with $0\ls \tau_i-\sigma_i\ls1$ for all $i=0,\ldots,k$. Let $W$ be an integer with
$$W\gs 1+\max\{\tau_i-\tau_{i+1}:\ i=0,\ldots,k-1\}.$$

{\rm (i)} For any integer $t$ with $t\ls\tau_0-1$ or $t\gs\tau_k+kW$, we have
\begin{equation}\label{t-ineq}\prod_{i=0}^k(t-\sigma_i-iW)(t+1-\tau_i-iW)\gs0.
\end{equation}

{\rm (ii)} Every interval $[\sigma_i,\tau_i]\ (0\ls i\ls k)$ contains an integer if and only if
\eqref{t-ineq} holds for all $t\in\Z$.
\end{lemma}
\Proof. Set $\sigma_i'=\sigma_i+iW$ and $\tau_i'=\tau_i+iW$ for all $i=0,\ldots,k$. Note that
$$\tau_i'-1\ls\sigma_i'\ls\tau_i'.$$
If $0\ls i<k$, then
$$\tau_i'=\tau_i+iW\ls\tau_{i+1}-1+W+iW=\tau_{i+1}'-1\ls\sigma_{i+1}'.$$

Let $t\in\Z$ and $0\ls j\ls k$. Suppose that $(t-\sigma_j')(t+1-\tau_j')<0$, which is equivalent to
$\tau_j'-1<t<\sigma_j'$. For $0\ls i<j$ we have $\sigma_i'\ls\tau_{i+1}'-1\ls\cdots\ls\tau_j'-1<t$.
If $j<i\ls k$, then $t<\sigma_j'\ls\tau_{j+1}'-1\ls\cdots\ls \tau_i'-1\ls\sigma_i'$. Thus $(t-\sigma_i')(t+1-\tau_i')>0$ for all $i=0,\ldots,k$ with $i\not=j$, and hence
$$\prod_{i=0}^k(t-\sigma_i')(t+1-\tau_i')<0.$$
Therefore,
\begin{align*}&\prod_{i=0}^k(t-\sigma_i')(t+1-\tau_i')\gs0
\\\iff& (t-\sigma_j')(t+1-\tau_j')\gs0 \ \t{for all}\ j=0,\ldots,k
\\\iff&t\not\in\bigcup_{j=0}^k(\tau_j'-1,\sigma_j').
\end{align*}

(i) If $t\ls\tau_0-1$, then $t\ls\tau_0'-1\ls\tau_1'-1\ls\cdots\ls\tau_k'-1$ and hence
$\prod_{i=0}^k(t-\sigma_i')(t+1-\tau_i')\gs0$ by the above. Similarly, if $t\gs\tau_k+kW$ then
$\sigma_0'\ls\sigma_1'\ls\cdots\ls\sigma_k'\ls\tau_k'\ls t$ and hence $\prod_{i=0}^k(t-\sigma_i')(t+1-\tau_i')\gs0$.

(ii) For any $i=0,\ldots,k$, clearly
$(\tau_i-1,\sigma_i)\cup[\sigma_i,\tau_i]=(\tau_i-1,\tau_i]$ contains a unique integer. So
\begin{align*} &[\sigma_i,\tau_i]\ \t{contains an integer for all}\ i=0,\ldots,k
\\\iff&(\tau_i-1,\sigma_i)\ \t{contains no integer for all}\ i=0,\ldots,k
\\\iff&(\tau_i'-1,\sigma_i')\ \t{contains no integer for all}\ i=0,\ldots,k
\\\iff&\prod_{i=0}^k(t-\sigma_i')(t+1-\tau_i')\gs0\quad\t{for all}\ t\in\Z.
\end{align*}

In view of the above, we have completed the proof of Lemma \ref{Lem2.2}.  \qed

\begin{lemma}\label{Lem2.3} Let $\da$ and $L$ be positive integers. Suppose that $z_0,\ldots,z_{\nu}\in\N$ and
$$P(z_0,\ldots,z_{\nu})=\sum_{i_0,\ldots,i_{\nu}\in\N\atop i_0+\cdots+i_{\nu}\ls \da}a_{i_0,\ldots,i_{\nu}}z_0^{i_0}\ldots z_{\nu}^{i_{\nu}}$$
with $a_{i_0,\ldots,i_{\nu}}\in\Z$ and $|a_{i_0,\ldots,i_{\nu}}|\ls L$.
Let $B$ be any integer greater than $2(1+z_0+\cdots+z_{\nu})^{\da}\da!L$.
Then $P(z_0,\ldots,z_{\nu})=0$ if and only if
\begin{equation}\label{CDz}\f{2C(B)D(B)-B^{(\da+1)^{\nu+1}}}{2B^{(\da+1)^{\nu+1}+1}}\ls z\ls \f{2C(B)D(B)+B^{(\da+1)^{\nu+1}}}{2B^{(\da+1)^{\nu+1}+1}}
\end{equation}
for some integer $z$, where
$C(x)=(1+\sum_{i=0}^{\nu}z_ix^{(\da+1)^i})^{\da}$
and
$$D(x)=\sum_{i_0,\ldots,i_{\nu}\in\N\atop i_0+\cdots+i_{\nu}\ls\da}i_0!\ldots i_{\nu}!(\da-i_0-\cdots-i_{\nu})!a_{i_0,\ldots,i_{\nu}}x^{(\da+1)^{\nu+1}-\sum_{j=0}^{\nu}i_j(\da+1)^j}.$$
\end{lemma}
\Proof. Write
$$C(x)=\sum_{i=0}^{\da(\da+1)^{\nu}}c_ix^i\ \t{and}\ D(x)=\sum_{j=0}^{(\da+1)^{\nu+1}}d_jx^j.$$
Then $c_i\gs0$, and also $|d_j|\ls\da!L$ since the multi-nomial coefficient
$$\bi{\da}{i_0,\ldots,i_{\nu},\da-i_0-\cdots-i_{\nu}}=\f{\da !}{i_0!\ldots i_{\nu}!(\da-i_0-\cdots-i_{\nu})!}
\gs1$$
for all $i_0,\ldots,i_{\nu}\in\N$ with $i_0+\cdots+i_{\nu}\ls\da$.
Write
$$C(x)D(x)=\sum_{k=0}^{(2\da+1)(\da+1)^{\nu}}r_kx^k.$$
Then
$$r_k=\sum_{0\ls i\ls \da(\da+1)^{\nu}\atop {0\ls j\ls(\da+1)^{\nu+1}\atop i+j=k}}c_id_j$$
and
$$ |r_k|\ls \sum_{i=0}^{\da(\da+1)^{\nu}}c_i\da!L=C(1)\da!L=(1+z_0+\cdots+z_{\nu})^{\da}\da!L<\f B2.$$
By the multi-nomial theorem,
$$C(x)=\sum_{i_0,\ldots,i_{\nu}\in\N\atop i_0+\cdots+i_{\nu}\ls\da}\bi{\da}{i_0,\ldots,i_{\nu},\da-i_0-\cdots-i_{\nu}}
z_0^{i_0}\ldots  z_{\nu}^{i_{\nu}}x^{\sum_{j=0}^{\nu}i_j(\da+1)^j}.$$
Therefore
$$r_{(\da+1)^{\nu+1}}=\sum_{i_0,\ldots,i_{\nu}\in\N\atop i_0+\cdots+i_{\nu}\ls\da}\da!a_{i_0,\ldots,i_{\nu}}z_0^{i_0}\ldots z_{\nu}^{i_{\nu}}=\da!P(z_0,\ldots,z_{\nu}).$$

Suppose that $P(z_0,\ldots,z_{\nu})=0$. Then $r_{(\da+1)^{\nu+1}}=0$. For the integer
$$z=\sum_{k=(\da+1)^{\nu+1}+1}^{(2\da+1)(\da+1)^{\nu}}r_kB^{k-1-(\da+1)^{\nu+1}},$$
we have
$$C(B)D(B)-zB^{(\da+1)^{\nu+1}+1}=\sum_{k=0}^{(\da+1)^{\nu+1}-1}r_kB^k$$
and hence
\begin{align*}\l|C(B)D(B)-zB^{(\da+1)^{\nu+1}+1}\r|\ls&\sum_{k=0}^{(\da+1)^{\nu+1}-1}|r_k|B^k
\\\ls&\f{B-1}2\sum_{k=0}^{(\da+1)^{\nu+1}-1}B^k\ls\f{B^{(\da+1)^{\nu+1}}}2.
\end{align*}
Therefore \eqref{CDz} is valid.

Now we assume that \eqref{CDz} holds for some $z\in\Z$. We want to show that $P(z_0,\ldots,z_{\nu})=0$.
By \eqref{CDz} we have
$$\l|C(B)D(B)-zB^{(\da+1)^{\nu+1}+1}\r|\ls\f12B^{(\da+1)^{\nu+1}}.$$
Let
$$S:=r_{(\da+1)^{\nu+1}}+\sum_{k=(\da+1)^{\nu+1}+1}^{(2\da+1)(\da+1)^{\nu}}r_kB^{k-(\da+1)^{\nu+1}}-Bz.$$
Then
\begin{align*}SB^{(\da+1)^{\nu+1}}=&r_{(\da+1)^{\nu+1}}B^{(\da+1)^{\nu+1}}
+\sum_{k=(\da+1)^{\nu+1}+1}^{(2\da+1)(\da+1)^{\nu}}r_kB^k-zB^{(\da+1)^{\nu+1}+1}
\\=&C(B)D(B)-\sum_{k=0}^{(\da+1)^{\nu+1}-1}r_kB^k-zB^{(\da+1)^{\nu+1}+1}
\end{align*}
and hence
\begin{align*}|SB^{(\da+1)^{\nu+1}}|\ls&|C(B)D(B)-zB^{(\da+1)^{\nu+1}+1}|+\sum_{k=0}^{(\da+1)^{\nu+1}-1}|r_k|B^k
\\\ls&\f{B^{(\da+1)^{\nu+1}}}2+\f{B-1}2\sum_{k=0}^{(\da+1)^{\nu+1}-1}B^k=B^{(\da+1)^{\nu+1}}-\f12.
\end{align*}
Therefore we must have $S=0$. It follows that $B\mid r_{(\da+1)^{\nu+1}}$. As
$$\l|r_{(\da+1)^{\nu+1}}\r|<\f B2,$$
we have
$$\da!P(z_0,\ldots,z_{\nu})=r_{(\da+1)^{\nu+1}}=0$$
and hence $P(z_0,\ldots,z_{\nu})=0$.

In view of the above, we have completed the proof of Lemma \ref{Lem2.3}. \qed

Lemma \ref{Lem2.3} and its proof also appeared in the author's recent book \cite[pp.\,117-119]{S24} in Chinese.

Now we present our auxiliary theorem.

\begin{theorem}\label{Th2.1} Let $A\se \N$ be any r.e. set. Then, there are $L[x]\in\Z[x]$
 and $M(x,y,z,t)\in\Z[x,y,z,t]$ satisfying the following {\rm (i)-(iii)}.

 {\rm (i)} $L(a)>0$ for all $a\in\Z$.

 {\rm (ii)} There are $k_0,k_1,k_2\in\Z^+$ such that $M(a,b,c,t)\gs0$
 whenever $a,b,c,t$ are integers with $a\gs0$, $b>1$, and $t<-c^2\lor t>R(a,c)$,
 where $R(a,c)=k_0(1+c)^{2k_1}L(a)+k_2$.

 {\rm (iii)} For any $a\in \N$ and any infinite subset $S$ of $\N$, we have
 $$a\in A\iff \exists b\in S\exists c\forall t\,(M(a,b,c,t)\gs0).$$

 {\rm (iv)} For any infinite subset $S$ of $\N$, there is a positive integer $n$ such that for any $a\in A$ and $N\in\N$
 there are $b\in S$ and $c\in\Z$ for which
 $b\gs N$, $b\mid c$, $0<c<b^n$, and $M(a,b,c,t)\gs0$ for all $t\in\Z$.
 \end{theorem}
 \Proof. By Matiyasevich's theorem \cite{M70}, there is a polynomial $P_0(z_0,z_1,\ldots,z_{\nu})\in\Z[z_0,\ldots,z_{\nu}]$ such that
 for any $a\in\N$ we have
 $$a\in A\iff \exists z_1\gs0\ldots\exists z_{\nu}\gs0\,(P_0(a,z_1,\ldots,z_{\nu})=0).$$
 Define
 $$P(a,z_0,\ldots,z_{\nu})=(z_0-1)^2+P_0^2(a,z_1,\ldots,z_{\nu}),$$
 and write
 $$P(a,z_0,\ldots,z_{\nu})=\sum_{i_0,\ldots,i_{\nu}\in\N\atop i_0+\cdots+i_{\nu}\ls \da}p_{i_0,\ldots,i_{\nu}}(a)z_0^{i_0}\ldots z_{\nu}^{i_{\nu}},$$
 where $\da$ is a positive even number. Set
 $$L(a)=\sum_{i_0,\ldots,i_{\nu}\in\N\atop i_0+\cdots+i_{\nu}\ls\da}p_{i_0,\ldots,i_{\nu}}(a)^2.$$
 Then $L(a)\gs p_{2,0,\ldots,0}(a)^2=1$ for all $a\in\Z$. As $\da$ is even, we have
 $$B(a,b):=2(\nu+1)^{\da}b^{\da}\da!L(a)\gs2$$
 for all $a,b\in\Z$ with $b\not=0$.

 Now fix $a\in\N$ and $b\in\{2,3,\ldots\}$. Then $B(a,b)\gs b^{\da}\gs b\gs2$.
 For convenience, we set $n_j=(\da+1)^j$ for all $j=0,\ldots,\nu+1$.
 Note that $P_0(a,z_1,\ldots,z_{\nu})=0$ for some $z_1,\ldots,z_{\nu}\in\{0,\ldots,b-1\}$ if and only if
$P(a,z_0,\ldots,z_{\nu})=0$ for some $z_0,\ldots,z_{\nu}\in\{0,\ldots,b-1\}$.
If $z_0,\ldots,z_{\nu}\in\{0,\ldots,b-1\}$, then
$$B(a,b)>2(1+(\nu+1)(b-1))^{\da}\da!L(a)\gs 2(1+z_0+\cdots+z_{\nu})^{\da}\da!L(a).$$
Thus, in view of Lemma \ref{Lem2.3},
$P(a,z_0,\ldots,z_{\nu})=0$ for some $z_0,\ldots,z_{\nu}\in\{0,\ldots,b-1\}$ if and only if for some
$$c\in\l\{\sum_{i=0}^{\nu}z_iB(a,b)^{n_i}:\ z_0,\ldots,z_{\nu}\in\{0,\ldots,b-1\}\r\}$$ we have
$$\f{2(1+c)^{\da}D(a,b)-B(a,b)^{n_{\nu+1}}}{2B(a,b)^{n_{\nu+1}+1}}\ls z\ls
\f{2(1+c)^{\da}D(a,b)+B(a,b)^{n_{\nu+1}}}{2B(a,b)^{n_{\nu+1}+1}}$$
for some integer $z$, where $$D(a,b)=\sum_{i_0,\ldots,i_{\nu}\in\N\atop i_0+\cdots+i_{\nu}\ls \da}i_0!\ldots i_{\nu}!(\da-i_0-\cdots-i_{\nu})!p_{i_0,\ldots,i_{\nu}}(a)B(a,b)^{n_{\nu+1}-\sum_{j=0}^{\nu}i_jn_j}.$$
Combining this with Lemma \ref{Lem2.1}, we see that
$P_0(a,z_1,\ldots,z_{\nu})=0$ for some $z_0,\ldots,z_{\nu}\in\{0,\ldots,b-1\}$ if and only if
for some $c\in\Z$, every interval $[\sigma_i,\tau_i]\ (i=0,\ldots,\nu+2)$ contains an integer, where
\begin{gather*}\sigma_0=\tau_0=\f c{B(a,b)},
\\\sigma_i=\f{c+1-bB(a,b)^{n_{i-1}}}{B(a,b)^{n_i}}\ \t{and}\
\tau_i=\f c{B(a,b)^{n_i}}\ (i=1,\ldots,\nu),
\\\sigma_{\nu+1}=\f{c+1-bB(a,b)^{n_{\nu}}}{(b^2+c^2)B(a,b)^{n_{\nu}}}
\ \t{and}\ \tau_{\nu+1}=\f c{(b^2+c^2)B(a,b)^{n_{\nu}}},
\\\sigma_{\nu+2}=\f{2(1+c)^\da D(a,b)-B(a,b)^{n_{\nu+1}}}{2B(a,b)^{n_{\nu+1}+1}}
\ \t{and}\ \tau_{\nu+2}=\f{2(1+c)^\da D(a,b)+B(a,b)^{n_{\nu+1}}}{2B(a,b)^{n_{\nu+1}+1}}.
\end{gather*}
Observe that
$$|\tau_i-\tau_{i+1}|=\l(\f1{B(a,b)^{n_i}}-\f1{B(a,b)^{n_{i+1}}}\r)|c|\ls\f{|c|}{B(a,b)^{n_i}}\ls\f{|c|}2$$
for all $i=0,\ldots,\nu-1$, and
$$|\tau_{\nu}-\tau_{\nu+1}|=\l(1-\f1{b^2+c^2}\r)\f{|c|}{B(a,b)^{n_{\nu}}}
\ls \f{|c|}{B(a,b)^{n_{\nu}}}\ls\f{|c|}2.$$
Note also that
$$|D(a,b)|\ls\sum_{i=0}^{n_{\nu+1}}\da!L(a)B(a,b)^i=\da!L(a)\f{B(a,b)^{n_{\nu+1}+1}-1}{B(a,b)-1}\ls \da!L(a)B(a,b)^{n_{\nu+1}+1}$$
and hence
\begin{align*}\tau_{\nu+1}-\tau_{\nu+2}\ls&\f{c/(b^2+c^2)}{B(a,b)^{n_{\nu}}}-\l(\f1{2B(a,b)}-\da!L(a)(1+c)^\da\r)
\\\ls&\f1{2bB(a,b)^{n_{\nu}}}-\f1{2B(a,b)}+\da!L(a)(1+c)^\da\ls1+\da!L(a)(1+c)^{\da}.
\end{align*}
Let $W=2+(1+c)^\da\da!L(a)$. Then $W-1\gs1+(1+c)^\da\gs 1+|c+1|\gs|c|\gs|c|/2$, and hence by the above we have
$$W\gs1+\max\{\tau_i-\tau_{i+1}:\ i=0,\ldots,\nu+1\}.$$

In view of the above and Lemma \ref{Lem2.2}(ii), $P_0(a,z_1,\ldots,z_{\nu})=0$ for some $z_1,\ldots,z_{\nu}
\in\{0,\ldots,b-1\}$ if and only if for some integer $c$ we have $Q(a,b,c,t)\gs0$ for all $t\in\Z$, where $Q(a,b,c,t)$ denotes
\begin{align*}&(B(a,b)t-c)(B(a,b)(t+1)-c)
\\\times&\prod_{i=1}^{\nu}\l(B(a,b)^{n_i}(t-iW)-c-1+bB(a,b)^{n_{i-1}}\r)\l(B(a,b)^{n_i}(t+1-iW)-c\r)
\\\times&\l((b^2+c^2)B(a,b)^{n_{\nu}}(t-(\nu+1)W)-c-1+bB(a,b)^{n_{\nu}}\r)
\\\times&\l((b^2+c^2)B(a,b)^{n_{\nu}}(t+1-(\nu+1)W)-c\r)
\\\times&\l(2B(a,b)^{n_{\nu+1}+1}(t-(\nu+2)W)-2(1+c)^{\da}D(a,b)+B(a,b)^{n_{\nu+1}}\r)
\\\times&\l(2B(a,b)^{n_{\nu+1}+1}(t+1-(\nu+2)W)-2(1+c)^{\da}D(a,b)-B(a,b)^{n_{\nu+1}}\r).
\end{align*}

Let $S$ be any infinite subset of $\N$. By the above, for any $a\in\N$ we have
\begin{align*}a\in A\iff& P_0(a,z_1,\ldots,z_{\nu})=0\ \t{for some}\ z_1,\ldots,z_{\nu}\in\N
\\\iff&\exists b\in S\,(b\gs2\land \exists z_1\in[0,b)\ldots\exists z_{\nu}\in[0,b)(P(a,z_1,\ldots,z_{\nu})=0))
\\\iff&\exists b\in S(b^2>b\land\exists c\forall t(Q(a,b,c,t)\gs0))
\\\iff&\exists b\in S\exists c\forall t\,(M(a,b,c,t)\gs0),
\end{align*}
where $M(a,b,c,t)=(b^2-b)(Q(a,b,c,t)+1)-1\in\Z[a,b,c,t]$ does not depend on $S$.

Given $a\in A$ and $N\in\N$, we may take $b\in S$ with
$$b\gs\max\{N,\,2(\nu+1)^{\da}\da!L(a)\}$$
such that $P(a,z_0,\ldots,z_{\nu})=0$ for some $z_0,\ldots,z_{\nu}\in[0,b)$ with $z_0=1$.
Then $c=\sum_{i=0}^{\nu}z_iB(a,b)^{n_{i}}\eq0\pmod b$ since $b\mid B(a,b)$.
By Lemmas 2.1-2.3, for any $t\in\Z$ we have
\begin{equation}
\label{t-ineq}\prod_{i=0}^{\nu+2}(t-\sigma_i-iW)(t+1-\tau_i-iW)\gs0\ \ (\t{i.e.,}\ Q(a,b,c,t)\gs0)
\end{equation}
It follows that  $M(a,b,c,t)\gs0$ for all $t\in\Z$. Note that
\begin{align*}1=z_0\ls c\ls&\sum_{i=0}^{\nu}(b-1)B(a,b)^{n_i}
\ls(b-1)\f{B(a,b)^{n_{\nu}+1}-1}{B(a,b)-1}
\\<&B(a,b)^{n_{\nu}+1}=(2(\nu+1)^{\da}b^\da \da!L(a))^{n_{\nu}+1}\ls (b^{\da+1})^{n_{\nu}+1}=b^n,
\end{align*}
where $n=(\da+1)(n_{\nu}+1)$ only depends on $A$.

Now it remains to show that (ii) in Theorem \ref{Th2.1} holds.
Let $a\in\N$, $b\in\{2,3,\ldots\}$ and $c\in\Z$. Then
$$-c^2-1\ls-|c|-1\ls\f{c}{B(a,b)}-1=\tau_0-1$$
and
\begin{align*}\tau_{\nu+2}+(\nu+2)W=&\f1{2B(a,b)}+\f{(1+c)^{\da}D(a,b)}{B(a,b)^{n_{\nu+1}+1}}+(\nu+2)W
\\\ls&\f1{2B(a,b)}+(1+c)^\da \da!L(a)+(\nu+2)W
\\<&1+(1+c)^\da \da!L(a)+(\nu+2)((1+c)^\da \da!L(a)+2)
\\=&(\nu+3)(1+c)^\da\da!L(a)+2\nu+5.
\end{align*} Thus, if $t$ is an integer with $t<-c^2$ or $t>R(a,c)=(\nu+3)(1+c)^\da \da!L(a)+2\nu+4$ then
by Lemma \ref{Lem2.2}(i) we have \eqref{t-ineq}
and hence $M(a,b,c,t)\gs0$. This concludes our proof. \qed

\section{Proof of Theorem \ref{Th1.1}}
 \setcounter{equation}{0}
 \setcounter{conjecture}{0}
 \setcounter{theorem}{0}

 For convenience, we define $\square=\{m^2:\ m\in\Z\}$.

 \begin{lemma}\label{Lem3.0} Let $C\in\Z$. Then
 \begin{align}\label{>=0}C\gs0&\iff\exists x\exists y\exists z(C=x^2+y^2+z^2+z),
 \\\label{Pell}C\gs0&\iff \exists_{x\not=0}((4C+2)x^2+1\in\square),
 \\\label{not=0}C\not=0&\iff \exists u\exists v(C=(2u+1)(3v+1)).
 \end{align}
 \end{lemma}
 \Proof. This is easy and known. Concerning \eqref{>=0}, by the Gauss-Legendre theorem on sums of three squares,
 $C\gs0$ if and only if $4C+1=(2x)^2+(2y)^2+(2z+1)^2$ (i.e., $C=x^2+y^2+z^2+z$) for some $x,y,z\in\Z$.
 By the theory of Pell equations, we have \eqref{Pell} which was first used by Sun \cite{S92b}.
 As any nonzero
 integer has the form $\pm3^a(3q+1)$ with $a\in\N$ and $q\in\Z$, we immediately get \eqref{not=0} which was
 an observation  due to Tung \cite{T85}. \qed

 \begin{lemma}\label{Lem3.1} Let $C_1,\ldots,C_n\in\Z$.

 {\rm (i)} We have
 \begin{align*}&C_1\gs0\lor\cdots\lor C_n\gs0
 \\\iff&\exists x\not=0((4C_1+2)x^2+1\in\square\lor\cdots\lor (4C_n+2)x^2+1\in\square)
 \\\iff&\exists u\exists v\exists w\l(\prod_{i=1}^n(2(2C_i+1)(2u+1)^2(3v+1)^2-w^2+1)=0\r).
 \end{align*}
 Also,
 \begin{align*}&C_1\gs0\land\cdots\land C_n\gs0
 \\\iff&\forall x\not=0\forall y\l(\prod_{i=1}^n((4C_i+2)x^2+y^2-1)\not=0\r)
 \\\iff&\forall x\forall y\exists u\exists v\l(x\l(\prod_{i=1}^n((4C_i+2)x^2+y^2-1)-(2u+1)(3v+1)\r)=0\r).
 \end{align*}

 {\rm (ii)} Suppose that $D_i\in\N$ and $|C_i|\ls D_i$ for all $i=1,\ldots,n$. Then
 \begin{align*}&C_1\gs0\land\cdots\land C_n\gs0
 \\\iff&\forall x\in[0,D_1\cdots D_n]\l(\prod_{i=1}^n(x+C_i+1)\not=0\r)
 \\\iff&\forall x\in[0,D_1\cdots D_n]\exists y\exists z\l(\prod_{i=1}^n(x+C_i+1)-(2y+1)(3z+1)=0\r)
 \end{align*}
 \end{lemma}
 \Proof. (i) The first assertion follows immediately from Lemma \ref{Lem3.0}.
 As for the second assertion, it suffices to note that
 $$C_i\gs0\iff -C_i-1\not\gs0\iff \forall x\not=0((4(-C_i-1)+2)x^2+1\not\in\square).$$

 (ii) If $C_i\gs 0$ for all $i=1,\ldots,n$, then for any $x\gs0$ we have $x+C_i+1>0$
 for all $i=1,\ldots,n$, and hence $\prod_{i=1}^n(x+C_i+1)\not=0$.
 If $C_i<0$ for some $1\ls i\ls n$, then for $x=-C_i-1$ we have $0\ls x\ls |C_i|\ls D_i\ls D_1\cdots D_n$. So part (ii) of Lemma \ref{Lem3.1} holds.

 In view of the above, we have completed the proof of Lemma \ref{Lem3.1}. \qed

\medskip
\noindent
{\it Proof of Theorem \ref{Th1.1}}. Take polynomials $R$ and $M$ as in Theorem \ref{Th2.1} and note that $S=\{b^2+2:\ b\in\N\}$ is an infinite set. By Theorem \ref{Th2.1}, for any $a\in\N$ we have
\begin{align*}a\in A\iff& \exists b\exists c\forall t[M(a,b^2+2,c,t)\gs0]
\\\iff&\exists b\exists c\forall t\in[-c^2,R(a,c)](M(a,b^2+2,c,t)\gs0)
\end{align*}
and hence
\begin{align*}a\in \bar A\iff& \forall b\forall c\exists t\,(-M(a,b^2+2,c,t)-1\gs0),
\end{align*}
also $M(a,b^2+2,c,t)\gs0$ whenever $t<-c^2$ or $t>R(a,c)$.
Moreover, if $a\in A$ then we may require further that $c>0$ and $(b^2+2)\mid c$.

(i) In view of the above and Lemmas \ref{Lem3.0}-\ref{Lem3.1}, we see that
$A$ has a $\exists^2\forall\exists^3$-representation over $\Z$ with $\forall$ bounded or unbounded,
and also a $\exists^2\forall^3\exists^2$-representation over $\Z$.
For $a\in\N$, $b,c\in\Z$ and $t\in [-c^2,R(a,c)]$, we clearly have $|M(a,b^2+2,c,t)|\ls P(a,b,c)^2$
for some $P(x,y,z)\in\Z[x,y,z]$. So, by using Lemma \ref{Lem3.1}(ii) we see that
$A$ also  has a $\exists^2\forall^2\exists^2$-representation over $\Z$ with $\forall$ bounded.
With the help of Lemma   \ref{Lem3.1}, $\bar A$ has a
$\forall^2\exists^4$-representation and also a $\forall^2\exists\forall^2\exists^2$-representation over $\Z$.

(ii) Let $D(c,s)=(s-c^2)(s-c^2-2c)$ and $a\in \N$. We claim that
\begin{align*}a\in A\iff& \exists s\forall t\exists c\gs0\,(D(c,s)\ls 0\land M(a,(s-c^2)^2+2,c,t)\gs0)
\\\iff&\exists s\forall t\in[-s^2,R(a,s)]\exists c\gs0
\\&\quad (D(c,s)\ls 0\land M(a,(s-c^2)^2+2,c,t)\gs0).
\end{align*}

Now we prove the claim. If $a\in A$, then for some $b\in\N$ and $c\in\Z^+$ with $(b^2+2)\mid c$ we have
$M(a,b^2+2,c,t)\gs0$ for all $t\in\Z$. As $0\ls b\ls b^2\ls c\ls 2c$,
for $s=b+c^2$ we have $c^2\ls s\ls c^2+2c$ and hence $D(c,s)\ls0$, and also
$$M(a,(s-c^2)^2+2,c,t)=M(a,b^2+2,c,t)\gs0$$
for all $t\in\Z$.

Now suppose that $s\in\Z$ and that for any $t\in[-s^2,R(a,s)]$ there is a number $c\in\N$
with $D(c,s)\ls 0$ and $M(a,(s-c^2)^2+2,c,t)\gs0$. Note that $c^2\ls s\ls c^2+2c<(c+1)^2$.
So $c=\lfloor\sqrt s\rfloor$ does not depend on $t$. Set $b=s-\lfloor\sqrt s\rfloor^2$.
Then $$M(a,b^2+2,c,t)=M(a,(s-c^2)^2+2,c,t)\gs0$$ for all $t\in[-s^2,R(a,s)]$.
If $t<-s^2$ then $t<-s\ls -c^2$ and hence $M(a,b^2+2,c,t)\gs0$.
If $t>R(a,s)$ then $t>R(a,c)$ (since $s\gs c^2\gs c\gs0$) and hence $M(a,b^2+2,c,t)\gs0$.
Therefore $M(a,b^2+2,c,t)\gs0$ for all $t\in\Z$, and hence $a\in A$.
 This concludes the proof of the claim.

 In view of the proved claim, for any $a\in\N$ we have
 \begin{align*}a\in A\iff&\exists s \forall t\exists c\,(c\gs0\land -D(c,s)\gs0\land M(a,(s-c^2)^2+2,c,t)\gs0)
 \\\iff&\exists s \forall t\in[-s^2,R(a,s)]\exists c\,(c\gs0\land -D(c,s)\gs0
 \\&\quad \land M(a,(s-c^2)^2+2,c,t)\gs0)
 \end{align*}
 and hence
 \begin{align*}a\in \bar A\iff&\forall s \exists t\forall c\,(-c-1\gs0\lor D(c,s)-1\gs0
 \\&\qquad \lor -M(a,(s-c^2)^2+2,c,t)-1\gs0)
 \end{align*}
 Combining this with Lemma \ref{Lem3.1}, we find that
 $A$ has a $\exists\forall\exists \forall^2\exists^2$-representation over $\Z$
 and a $\exists\forall\exists\forall\exists^2$-representation over $\Z$ with $\forall$ bounded.
 Also, $\bar A$ has a $\forall\exists\forall\exists^3$-representation over $\Z$.

 In view of the above, we have completed the proof of Theorem \ref{Th1.1}. \qed

\section{Proof of Theorem \ref{Th1.2}}
 \setcounter{equation}{0}
 \setcounter{conjecture}{0}
 \setcounter{theorem}{0}

 Let $J_k(x_1,\ldots,x_k,x)$ be the polynomial
$$\prod_{\ve_1,\ldots,\ve_k\in\{\pm1\}}\l(x+\ve_1\sqrt{x_1}+\ve_2\sqrt{x_2}X+\cdots+\ve_k\sqrt{x_k}X^{k-1}\r)$$
with $X=1+\sum_{j=1}^k x_j^2$. This polynomial (in $x_1,\ldots,x_n,x$) with integer coefficients was introduced by
Matiyasevich and Robinson \cite{MR}. For fixed $A_1,\ldots,A_k\in\Z$, the monic polynomial $J_k(A_1,\ldots,A_k,x)$ is of degree $2^k$ in $x$.

 \begin{lemma} \label{Lem4.1}
 Let $A_1,\ldots,A_k\in\Z$.

 {\rm (i)} We have
 $$A_1,\ldots,A_k\in\square\iff \exists x(J_k(A_1,\ldots,A_k,x)=0).$$

 {\rm (ii)} If $S,T\in\Z$ and $S\not=0$, then
 \begin{align*}&A_1\in\square\land \cdots\land A_k\in\square\land S\mid T
 \\\iff&\exists x\l(S^{2^k}J_k\l(A_1,\ldots,A_k,x+\f{T}{S}\r)=0\r).
 \end{align*}

 {\rm (iii) (Matiyasevich-Robinson Relation-Combining Theorem \cite{MR})} If $R,S,T\in\Z$ and $S\not=0$, then
 \begin{align*}&A_1\in\square\land \cdots\land A_k\in\square\land S\mid T\land R>0
 \\\iff&\exists n\gs0\l((S^2(1-2R))^{2^k}J_k\l(A_1,\ldots,A_k,T^2+W^k+\f{S^2n+T^2}{S^2(1-2R)}\r)=0\r),
 \end{align*}
 where $W=1+\sum_{i=1}^kA_i^2$.
 \end{lemma}
 \begin{remark} Parts (i) and (iii) of Lemma \ref{Lem4.1}
 were due to Matiyasevich and Robinson \cite[Theorems 1-3]{MR}.
 Part (ii) was stated explicitly in \cite[Lemma 17]{S92a}; in fact, if $x_0+T/S$ is a rational zero of the monic polynomial
 $J_k(A_1,\ldots,A_k,x)$ then it is an integer since any rational algebraic integer must belong to $\Z$.
 \end{remark}

 \medskip
 \noindent{\it Proof of Theorem \ref{Th1.2}}. Take polynomials $R$ and $M$ as in Theorem \ref{Th2.1} and note that $S=\{4b+3:\ b\in\square\}$ is an infinite set. By Theorem \ref{Th2.1}, for any $a\in\N$ we have
 $$a\gs0\land b\in\square\land(t<-(c+1)^2\lor t>R(a,c+1))\Rightarrow M(a,4b+3,c+1,t)\gs0$$
and
\begin{align*}a\in A\iff& \exists b\in\square\,\exists c\forall t(M(a,4b+3,c+1,t)\gs0).
\end{align*}
Moreover, if $a\in A$ then we may choose $b\in\square$ and $c\gs0$ with $(4b+3)\mid (c+1)$
such that $M(a,4b+3,c+1,t)\gs0$ for all $t\in\Z$.

Let $a\in\N$. We claim that
\begin{align*}&a\in A
\\\iff&\exists s\forall t\exists c\,(s-c^2\in\square\land (4(s-c^2)+3)\mid (c+1)
\\&\qquad \land (c+1)^2(M(a,4(s-c^2)+3,c+1,t)+1)>0)
\\\iff&\exists s\forall t\in[-(s+1)^2,R(a,s+1)]\exists c\,(s-c^2\in\square\land (4(s-c^2)+3)\mid (c+1)
\\&\qquad \land (c+1)^2(M(a,4(s-c^2)+3,c+1,t)+1)>0).
\end{align*}
When $a\in A$, we may choose $b\in\square$ and $c\gs0$ for which
$$4b+3\mid c+1\land \forall t(M(a,4b+3,c+1,t)\gs0).$$
Take $s=b+c^2$, Then $s-c^2=b\in\square$, $4(s-c^2)+3=4b+3$ divides $c+1$, and
$$M(a,4(s-c^2)+3,c+1,t)=M(a,4b+3,c+1,t)\gs0$$
for all $t\in\Z$. Note that $(c+1)^2(M(a,4(s-c^2)+3,c+1,t)+1)>0$.

Now we prove the remaining direction of the claim. Suppose that $s\in\Z$ and that for any
$t\in[-(s+1)^2,R(a,s+1)]$ there is an integer $c(t)$ for which
\begin{gather*} s-c(t)^2\in\square,\ (4(s-c(t)^2)+3)\mid (c(t)+1),
\\(c(t)+1)^2(M(a,4(s-c(t)^2)+3,c(t)+1,t)+1)>0.
\end{gather*}
Clearly, $c(t)+1\not=0$ and
\begin{align*}c(t)^2\ls&\ s=(s-c(t)^2)+c(t)^2<4(s-c(t)^2)+3+c(t)^2
\\\ls&\ |c(t)+1|+c(t)^2\ls(|c(t)|+1)^2.
\end{align*}
Hence $s\gs0$ and $|c(t)|=\lfloor\sqrt s\rfloor$. Since
$$\lfloor \sqrt{s}\rfloor+1+(-\lfloor\sqrt s\rfloor+1)=2\not\eq0\pmod{4(s-c(t)^2)+3},$$
there is a unique $c\in\{\pm\lfloor \sqrt s\rfloor\}$ with $c+1$ divisible by
$4(s-\lfloor s\rfloor^2)+3$. It follows that $c(t)=c$ for all $t\in\Z$. Set $b=s-c^2=s-\lfloor \sqrt s\rfloor^2$. Then $b\in\square$, $4b+3\mid c+1$, and $M(a,4b+3,c+1,t)\gs0$ for all $t\in[-(s+1)^2,R(a,s+1)]$. If $t<-(s+1)^2$, then $t<-(\lfloor s\rfloor+1)^2\ls-(c+1)^2$
and hence $M(a,4b+3,c+1,t)\gs0$. Note that
$$(1+(\lfloor\sqrt s\rfloor+1))^2\gs (1-\lfloor \sqrt s\rfloor+1)^2.$$
If $t>R(a,s+1)$, then $t>R(a,\lfloor s\rfloor +1)\gs R(a,c+1)$ and hence
$M(a,4b+3,c+1,t)\gs0$.  As $b\in\square$ and $M(a,4b+3,c+1,t)\gs0$ for all $t\in\Z$,
we have $a\in A$. This concludes the proof of the claim.

Combining the proved claim with Lemma \ref{Lem4.1}(iii) and \eqref{>=0}, we get that
$A$ has a $\exists\forall \exists^4$-representation with $\forall$ bounded (or unbounded) over $\Z$.

By the proved claim, \eqref{Pell} and Lemma \ref{Lem4.1}(ii), for any $a\gs0$ we have
\begin{align*}&a\in A
\\\iff&\exists s\forall t\exists c(s-c^2\in\square\land (4(s-c^2)+3)\mid (c+1)
\\&\land \exists d\not=0((4(c+1)^2(M(a,4(s-c^2)+3,c+1,t)+1)-2)d^2+1\in\square)
\\\iff& \exists s\forall t\exists c\exists d\not=0\exists x(P(a,s,t,c,d,x)=0),
\end{align*}
where $P$ is a suitable polynomial with integer coefficients. It follows that
\begin{align*}a\in \bar A\iff& \forall s\exists t\forall c\forall d\not=0\forall x
(P(a,s,t,c,d,x)\ne0)
\\\iff&\forall s\exists t\forall c\forall d\forall x\exists y\exists z(d(P(a,s,t,c,d,x)-(2y+1)(3z+1))=0)
\end{align*}
with the aid of \eqref{not=0}. So $\bar A$ has a $\forall \exists\forall^3\exists^2$-representation over $\Z$. This concludes our proof of Theorem \ref{Th1.2}. \qed

\section{Proof of Theorem \ref{Th1.3}}
 \setcounter{equation}{0}
 \setcounter{conjecture}{0}
 \setcounter{theorem}{0}

 \begin{lemma}[Sun \cite{S92b}] \label{n+2} There is a polynomial $P(x_1,\ldots,x_{2n+2})$
 with integer coefficients
 such that for any $C_1,\ldots,C_n\in\Z$ we have
 \begin{align*}&C_1\gs0\land \cdots\land C_n\gs0
 \\\iff&\exists x_1\cdots\exists x_{n+2}(P(C_1,\ldots,C_n,x_1,\ldots,x_{n+2})=0).
 \end{align*}
 \end{lemma}

 \begin{lemma}\label{Lem5.1} There are polynomials
 $$P(x_1,\ldots,x_{2n+3})\in\Z[x_1,\ldots,x_{2n+3}]\ \t{and}\ Q(x_1,\ldots,x_{2n+2})\in\Z[x_1,\ldots,x_{2n+2}]$$
 such that for any $C_1,\ldots,C_n\in\Z$ we have
 \begin{align*}&C_1\gs0\lor\ldots\lor C_n\gs0
 \\\iff&\forall x_1\ldots\forall x_n\forall x\exists y\exists z(P(C_1,\ldots,C_n,x_1,\ldots,x_n,x,y,z)=0),
 \end{align*}
 and
 \begin{align*}&C_1\gs0\lor\ldots\lor C_n\gs0
 \\\iff&\forall x_1\in[0,D_1]\ldots\forall x_n\in[0,D_n]\exists y\exists z(P(C_1,\ldots,C_n,x_1,\ldots,x_n,y,z)=0)
 \end{align*}
 provided that $|C_i|\ls D_i$ with $D_i\in\N$ for all $i=1,\ldots,n$.
 \end{lemma}
 \Proof. (i) For each $i=1,\ldots, n$, clearly
 $$C_i<0\iff -C_i-1\gs0\iff \exists x_i\not=0(1-(4C_i+2)x_i^2\in\square).$$
 Thus
 \begin{align*}&\neg(C_1\gs0\lor\cdots\lor C_n\gs0)
 \\\iff&C_1<0\land\cdots\land C_n<0
 \\\iff&\exists x_1\not=0(1-(4C_1+2)x_1^2\in\square)\land\cdots\land \exists x_n\not=0(1-(4C_n+2)x_1^2\in\square)
 \\\iff&\exists x_1\cdots\exists x_n(x_1\cdots x_n\not=0\land (1-(4C_1+2)x_1^2\in\square)
 \\&\quad \land\cdots\land (1-(4C_n+2)x_n^2\in\square))
 \\\iff&\exists x_1\cdots\exists x_n(x_1\cdots x_n\not=0
 \\&\qquad \land
 \exists x(J_n(1-(4C_1+2)x_1^2,\ldots,1-(4C_n+2)x_n^2,x)=0)
 \end{align*}
 and hence
 \begin{align*}&C_1\gs0\lor\cdots\lor C_n\gs0
 \\\iff&\forall x_1\cdots\forall x_n\forall x\,(x_1\cdots x_n=0
 \\&\ \lor
 J_n(1-(4C_1+2)x_1^2,\ldots,1-(4C_n+2)x_n^2,x)\not=0)
 \\\iff&\forall x_1\cdots\forall x_n\forall x\exists y\exists z
 \,(x_1\cdots x_n
 \\&\times(J_n(1-(4C_1+2)x_1^2,\ldots,1-(4C_n+2)x_n^2,x) -(2y+1)(3z+1))=0).
 \end{align*}

 (ii) We now prove the latter assertion in Lemma \ref{Lem5.1}. Let $D_1,\ldots,D_n\in\N$ with $D_i\gs|C_i|$
 for all $i=1,\ldots,n$. By Lemma \ref{Lem3.1},
 $$C_i\gs0\iff \forall x_i\in[0,D_i](x_i+C_i+1\not=0).$$
 Thus
 \begin{align*}&C_1\gs0\lor\ldots\lor C_n\gs0
 \\\iff&\forall x_1\in[0,D_1]\ldots\forall x_n\in[0,D_n](x_1+C_1+1\not=0\lor\cdots\lor x_n+C_n+1\not=0)
 \\\iff&\forall x_1\in[0,D_1]\ldots\forall x_n\in[0,D_n]\exists y\exists z
 \\&\quad ((x_1+C_1+1)^2+\cdots+(x_n+C_n+1)^2=(2y+1)(3z+1)).
 \end{align*}
This ends the proof. \qed

 \begin{lemma}\label{Lem5.2} Let $k,m\in\Z$ with $k>0$, $2\mid k$ and $m\eq3\pmod4$. Then there is a unique $b\in\N$ such that $|m-b^k|=\min_{x\in\Z}|m-x^k|$. Moreover, for $b\in\Z$ we have
 $$|m-b^k|=\min_{x\in\Z}|m-x^k|\iff |m-b^k|<|m-(b\pm1)^k|.$$
 \end{lemma}
 \Proof. If $a,b\in\N$ and $|m-a^k|=|m-b^k|$ but $a\not=b$, then $m-a^k=-(m-b^k)$ and hence $2m=a^k+b^k$, thus $a\eq b\pmod 2$ and we get a contradiction since $2m$ is neither divisible by $4$
 nor congruent to $2$ modulo $8$. (Note that an odd square is congruent to $1$ modulo $8$.)
 So, there is a unique $b\in\N$ with $|m-b^k|=\min_{x\in\Z}|m-x^k|$.

 If $b\in\Z$ and $|m-b^k|=\min_{x\in\Z}|m-x^k|$, then $|m-b^k|<|m-(b\pm1)^k|$ as $|b\pm1|\not=|b|$.

Suppose that $b\in\Z$ and $|m-b^k|<|m-(b\pm1)^k|$. If $b=0$, then $|m|\ls|m-1|$, hence $m\ls0$ and
$$\min_{x\in\Z}|m-x^k|=\min_{x\in\Z}|-|m|-|x|^k|=|m|=|m-b^k|.$$
Now assume that $b\not=0$. Then $|b|\pm1\gs0$. Note that
$$|m-|b|^k|=|m-b^k|\ls|m-(|b|\pm1)^k|.$$
If $m\ls (|b|-1)^k$, then $m-|b|^k<m-(|b|-1)^k\ls0$ and hence $|m-|b|^k|>|m-(|b|-1)^k|$
which leads to a contradiction. If $m\gs(|b|+1)^k$, then $m-|b|^k>m-(|b|+1)^k\gs0$, which also leads
to a contradiction. Therefore
$$(|b|-1)^k<m<(|b|+1)^k.$$
If $x\in\Z$ and $|x|=|b|$, then $|m-x^k|=|m-b^k|$ since $k$ is even.
For $x\in\Z$ with $|x|<|b|$, we have
$m-x^k=m-|x|^k\gs m-(|b|-1)^k>0$ and hence $$|m-x^k|\gs|m-(|b|-1)^k|\gs|m-b^k|.$$
For $x\in\Z$ with $|x|>|b|$, we have
$m-x^k=m-|x|^k\ls m-(|b|+1)^k<0$ and hence $$|m-x^k|\gs|m-(|b|+1)^k|\gs|m-b^k|.$$
So we have $|m-b^k|=\min_{x\in\Z}|m-x^k|$.

The proof of Lemma \ref{Lem5.2} is now complete. \qed

\medskip
\noindent{\it Proof of Theorem \ref{Th1.3}}. Take polynomials $R$ and $M$ as in Theorem \ref{Th2.1} and note that $S=\{(2x)^2+4:\ x\in\Z\}$ is an infinite subset of $\N$. By Theorem \ref{Th2.1}, for any $a\in\N$ we have
 $$ b\in\Z\land c\in\Z\land (t<-c^2\lor t>R(a,c)\Rightarrow M(a,b^2+4,c,t)\gs0$$
and
\begin{align*}a\in A\iff& \exists b\exists c\forall t(M(a,b^2+4,c,t)\gs0)
\\\iff&\exists b\,(2\mid b\land \exists c\forall t\,(M(a,b^2+4,c,t)\gs0)).
\end{align*}
Moreover, if $a\in A$ then we may choose $b\gs2$ and $0<c<(b^2+4)^n$ with $(b^2+4)\mid c$
such that $M(a,b^2+4,c,t)\gs0$ for all $t\in\Z$, where $n$ is a positive integer only depending on $A$.

Note that $k=4n$ is a positive even number. For $b,q\in\Z$ let
$$P^+(b,q)=(4q-1-(b+1)^k)^2-(4q-1-b^k)^2$$
and
$$P^-(b,q)=(4q-1-(b-1)^k)^2-(4q-1-b^k)^2.$$
By Lemma \ref{Lem5.1},
$$|4q-1-b^k|=\min_{x\in\Z}|4q-1-x^k|\iff P^+(b,q)>0\land P^-(b,q)>0.$$

Let $a\in\N$. We claim that
\begin{align*}a\in A\iff&\exists q\forall b\forall t(P^+(b,q)\ls0\lor P^-(b,q)\ls0\lor M(a,b^2+4,4q-b^k,t)\gs0)
\\\iff&\exists q\forall b\in[0,8q^2+1]\forall t\in[-((4q-1)^2+1)^2,R(a,(4q-1)^2+1)]
\\&(P^+(b,q)\ls0\lor P^-(b,q)\ls0\lor M(a,b^2+4,4q-b^k,t)\gs0).
\end{align*}

Suppose that $a\in A$. Then there are $b_0,c\in\Z$ with
$$2\mid b_0,\ b_0^2+4\gs6,\ (b_0^2+4)\mid c\ \t{and}\ 0<c<(b_0^2+4)^n$$
such that $M(a,b_0^2+4,c,t)\gs0$ for all $t\in\Z$. As $4\mid b_0^2$ and $4\mid c$,
$q=(b_0^k+c)/4$ is an integer. Let $m=4q-1$. Note that
\begin{gather*}0\ls c-1<(b_0^2+4)^n\ls(2b_0^2)^n\ls|b_0|^{3n}\ls|b_0|^{4n-1},
\\2(c-1)<4n|b_0|^{4n-1}\ls (|b_0|+1)^{4n}-|b_0|^{4n},
\\|m-b_0^{4n}|=c-1<(|b_0|+1)^{4n}-(|b_0|^{4n}+c-1)=-(m-(|b_0|+1)^{4n}),
\\|m-b_0^{4n}|=c-1<c-1+|b_0|^{4n}-(|b_0|-1)^{4n}=m-(|b_0|-1)^{4n}.
\end{gather*}
Therefore $|m-b_0^k|<|m-(b_0\pm1)^k|$, hence $P^+(b_0,q)>0$ and $P^-(b_0,q)>0$.
If $b\in\Z$ and $P^{\pm}(b,q)>0$, then $|m-b^k|=\min_{x\in\Z}|m-x^k|=|m-b_0^k|$
and hence $|b|=|b_0|$, thus $4q-b^k=b_0^k+c-b^k=c$ and
$$M(a,b^2+4,4q-b^k,t)=M(a,b^2+4,c,t)\gs0$$
for all $t\in\Z$. So, for any $b\in\Z$ we have
$$P^+(b,q)\ls0\lor P^-(b,q)\ls 0\lor\forall t(M(a,b^2+4,4q-b^k,t)\gs0).$$

Now  we prove another direction of the claim. Let $q\in\Z$ and assume that for any $b\in[0,8q^2+1]$
and $t\in[-((4q-1)^2+1)^2,R(a,(4q-1)^2+1)]$ we have
$$P^+(b,q)\ls0\lor P^-(b,q)\ls0\lor M(a,b^2+4,4q-b^k,t)\gs0.$$
Take the unique $b\in\N$ with $|m-b^k|=\min_{x\in\Z}|m-x^k|$, where $m=4q-1$.
Then $|m-b^k|<|m-(b\pm1)^k|$ and hence both $P^+(b,q)$ and $P^-(b,q)$ are positive.
If $b\not=0$, then $b^k-|m|\ls|b^k-m|<|0^k-m|=|m|$. No matter $b=0$ or not, we have $b^k\ls 2|m|-1$.
Hence
 $$0\ls b\ls 2|m|-1\ls2(4|q|+1)-1\ls 8q^2+1.$$
 If $t\in[-(m^2+1)^2,R(a,m^2+1)]$, then by the assumption we have
 $M(a,b^2+4,c,t)\gs0$, where $c=4q-b^k$. Note that
 $$|c|=|m+1-b^k|\ls |m-b^k|+1\ls|m-0^k|+1\ls m^2+1$$
 and hence $|1+c|\ls1+|c|\ls 1+(m^2+1)$. If $t<-(m^2+1)^2$ or $t>R(a,m^2+1)$, then
 $t<-c^2$ or $t>R(a,c)$, and hence $M(a,b^2+4,c,t)\gs0$. So $M(a,b^2+4,c,t)\gs0$ for all $t\in\Z$,
 and hence $a\in A$. This concludes the proof of the claim.

By the proved claim, for any $a\in\N$ we have
\begin{align*}&a\in A
\\\iff&\exists q\forall b\forall t\,(-P^+(b,q)\gs0\lor -P^-(b,q)\gs0\lor M(a,b^2+4,4q-b^k,t)\gs0)
\\\iff&\exists q\forall b\in[0,8q^2+1]\forall t\in[-((4q-1)^2+1)^2,R(a,(4q-1)^2+1)]
\\&(-P^+(b,q)\gs0\lor -P^-(b,q)\gs0\lor M(a,b^2+4,4q-b^k,t)\gs0).
\end{align*}
Clearly $|P^{\pm}(b,q)|\ls P_0(q)^2$ for all $b\in[0,8q^2+1]$, and
$$|M(a,b^2+4,4q-b^k,t)|\ls M_0(a,q)^2$$ for all $b\in[0,8q^2+1]$ and
$t\in[-((4q-1)^2+1)^2,R(a,(4q-1)^2+1)]$,
where $P_0$ and $M_0$ are suitable polynomials with integer coefficients.

Combining the last paragraph with Lemma \ref{Lem3.1}(i), we find that $A$ has a $\exists\forall^2\exists^3$-representation over $\Z$ and also a $\exists\forall^2\exists^3$-representation over $\Z$ with $\forall$ bounded.
Combining the last paragraph with Lemma \ref{Lem5.1}, we see that that $A$ has a
$\exists \forall^6\exists^2$-representation over $\Z$, and also
a $\exists\forall^5\exists^2$-representation over $\Z$ with $\forall$ bounded.

Note that
\begin{align*}a\in \bar A=\N\sm A\iff&\forall q\exists b\exists t\,(P^+(b,q)-1\gs0\land P^-(b,q)-1\gs0
\\&\qquad\ \ \quad\land -M(a,b^2+4,4q-b^k,t)-1\gs0).
\end{align*}
Combining this with Lemma \ref{Lem3.1}(i), we get that $\bar A$ has a $\forall\exists^2\forall^2\exists^2$-representation over $\Z$;
 if we apply Lemma \ref{n+2},
 then we find that $\bar A$ has a $\forall\exists^7$-representation over $\Z$.

 The proof of Theorem \ref{Th1.3} is now complete. \qed

\Ack. The author would like to thank the referee for many helpful suggestions.


 \end{document}